\documentclass[12pt,a4paper]{article}
\usepackage{amsfonts}
\usepackage{mathrsfs}
\usepackage{hyperref}
\usepackage{xspace,colortbl}
\usepackage{amssymb}
\usepackage{amsmath}

\newtheorem{theorem}{Theorem}[section]

\newtheorem{lemma}[theorem]{Lemma}

\newtheorem{proposition}[theorem]{Proposition}

\newtheorem{remark}[theorem]{Remark}
\newenvironment{proof}{{\bf Proof.  }}{$\square$}

\begin{document}
\title{Gr\"{o}bner-Shirshov bases for free inverse semigroups\footnote{Supported by the
NNSF of China  (No.10771077) and the NSF of Guangdong Province
 (No.06025062).}}
\author{
L. A. Bokut\footnote {Supported by the grant LSS--344.2008.1 and SB
RAS
Integration grant No. 2006.1.9 (Russia).} \\
{\small \ School of Mathematical Sciences, South China Normal
University}\\
{\small Guangzhou 510631, P. R. China}\\
{\small Sobolev Institute of Mathematics, Russian Academy of
Sciences}\\
{\small Siberian Branch, Novosibirsk 630090, Russia}\\
{\small Email: bokut@math.nsc.ru}\\
\\
 Yuqun
Chen\footnote {Corresponding author.} \  and Xiangui Zhao\\
{\small \ School of Mathematical Sciences, South China Normal
University}\\
{\small Guangzhou 510631, P. R. China}\\
{\small Email: yqchen@scnu.edu.cn}\\
{\small xian.zhao@umanitoba.ca}}

\date{}

\maketitle \noindent\textbf{Abstract:}  A new construction for free
inverse semigroups was obtained by Poliakova and Schein in 2005.
Based on their result, we find Gr\"{o}bner-Shirshov bases for  free
inverse semigroups with respect to the deg-lex order of words. In
particular, we give the (unique and shortest) normal forms in the
classes of equivalent words of a free inverse semigroup together
with the Gr\"{o}bner-Shirshov algorithm to transform any word to its
normal form.



\section{Introduction}
The theories of Gr\"obner and Gr\"obner--Shirshov bases were
invented independently by A. I. Shirshov  \cite{Sh} for
non-commutative and non-associative algebras, and by H. Hironaka
\cite{H} and B. Buchberger \cite{bu65} for commutative algebras. The
technique of Gr\"obner--Shirshov bases is proved to be very useful
in the study of presentations of associative algebras, Lie algebras,
semigroups, groups,  $\Omega$-algebras, etc. by generators and
defining relations, see, for example, the book \cite{BK} by L. A.
Bokut and G. Kukin, survey papers \cite{BK00, BK03} by L. A. Bokut
and P. Kolesnikov, and \cite{BS} by L. A. Bokut and Y. Q. Chen.

Let us mention some recent results on Gr\"obner--Shirshov bases for
groups and semigroups. Gr\"obner--Shirshov bases for braid groups in
different sets of generators were found in \cite{bokut08-1},
\cite{bokut08-2}, \cite{bokut-cha-sh}  and \cite{bokut-fo-ke-sh}. In
particular, Artin-Markov (\cite{artin}, \cite{markov}), Garside
(\cite{garside}), and Birman-Ko-Lee (\cite{birman-ko-lee}) normal
forms of a braid group were given in these papers.
Gr\"obner--Shirshov basis for a Chinese monoid (\cite{duch-krob})
was found and the staircase normal form (\cite{duch-krob}) of a
Chinese monoid was given in \cite{chen-qiu}. Gr\"obner--Shirshov
basis for the Adjan extension of the Novikov group was found in
\cite{bokut-cha}, which also gave the Adjan proof of the Adjan-Rabin
theorem (\cite{adjan}, \cite{rabin}). Kalorkoti (\cite{kal82},
\cite{kal06}) has actually found Gr\"obner--Shirshov bases for some
groups of Novikov-Boone type (\cite{novikov}, \cite{boone}) and
given a new proof of Bokut  (\cite{bokut68}) and Collins'
(\cite{collins69}, \cite{collins72}) results.

Inverse semigroups form one of the important classes of semigroups.
If one treats them as unary semigroups with the additional unary
operation, taking the inverse element, then, as is known, the class
of inverse semigroups is a variety. Free inverse semigroups are the
free algebras of this variety. There are several works devoted to
different constructions describing free inverse semigroups, see H.
E. Scheiblich
 \cite{Scheib, Scheib73}, W. D. Munn  \cite{M, M74}, G. B. Preston
 \cite{Pr1973} and B. M. Schein  \cite{Schein1975}, see also a survey
 \cite{Re} by N. R. Reilly and a book \cite{Pe} by M. Petrich. A new
construction for free inverse semigroups was found recently in a
fundamental paper by O. Poliakova and B. M. Schein \cite{Schein}. As
was noted in \cite{Schein}, each of the constructions mentioned
above can be easily obtained by using the construction of Poliakova
and Schein.

In this paper we find  Gr\"obner--Shirshov bases for free inverse
semigroups using the concept of the canonical idempotents from the
paper \cite{Schein}. As a result, we obtain the (unique and
shortest) normal forms of elements of the free inverse semigroup
together with the Gr\"{o}bner-Shirshov algorithm to transform any
word to its normal form. The normal forms consists of a set of
canonical words in the sense of \cite{Schein}, but contrary to
canonical word, the normal form is unique for a given element of a
free inverse semigroup.

\section{Preliminaries}
We first cite some concepts and results from the literature
\cite{Sh, B72, B76} which are related to the Gr\"{o}bner-Shirshov
bases for associative algebras.

Let $k$ be a field, $k\langle X\rangle$ the free associative algebra
over $k$ generated by $X$ and $ X^{*}$ the free monoid generated by
$X$, where the empty word is the identity which is denoted by 1. For
a word $w\in X^*$, we denote the length (degree) of $w$ by $|w|$.
Let $X^*$ be a well ordered set and $f\in k\langle X\rangle$. Then
by $\bar{f}$ we denote the maximum monomial in $f$, which is also
called the {\it leading word} of $f$. We call $f$ {\it monic} if
$\bar{f}$ has coefficient 1.

A well ordering $<$ on $X^*$ is called {\it monomial} if it is
compatible with the multiplication of words, that is, for $u, v\in
X^*$, we have
$$
u < v \Rightarrow w_{1}uw_{2} < w_{1}vw_{2},  \ for \  all \
 w_{1}, \ w_{2}\in  X^*.
$$
A standard example of monomial ordering on $X^*$ is the {\it deg-lex
ordering} to compare two words first by degree and then
lexicographically, where $X$ is a well ordered set.

Let $f$ and $g$ be two monic polynomials in \textmd{k}$\langle
X\rangle$ and $<$ a monomial ordering on $X^*$. Then, there are two
kinds of {\it compositions}:

$ (i)$ If \ $w$ is a word such that $w=\bar{f}b=a\bar{g}$ for some
$a,b\in X^*$ with $|\bar{f}|+|\bar{g}|>|w|$, then the polynomial
 $ (f,g)_w=fb-ag$ is called the {\it intersection composition} of $f$ and
$g$ with respect to $w$.

$ (ii)$ If  $w=\bar{f}=a\bar{g}b$ for some $a,b\in X^*$, then the
polynomial $ (f,g)_w=f - agb$ is called the {\it inclusion
composition} of $f$ and $g$ with respect to $w$.

The word $w$ is called the \emph{ambiguity} of the composition
$(f,g)_w$.

Let $S\subset k\langle X\rangle$ such that every $s\in S$ is monic.
Then the composition $ (f,g)_w$ is called {\it trivial} modulo
$(S,w)$  if $ (f,g)_w=\sum\alpha_i a_i s_i b_i$, where each
$\alpha_i\in k$, $a_i,b_i\in X^{*}, \ s_i\in S$ and $a_i
\overline{s_i} b_i<w$. If this is the case, then we write
$$
 (f,g)_w\equiv0\quad mod (S,w).
$$
In general, for $p,q\in k\langle X\rangle$, we write $ p\equiv
q\quad mod (S,w) $ which means that $p-q=\sum\alpha_i a_i s_i b_i $,
where each $\alpha_i\in k,a_i,b_i\in X^{*}, \ s_i\in S$ and $a_i
\overline{s_i} b_i<w$.

$S$ is called a {\it Gr\"{o}bner-Shirshov basis} in $k\langle
X\rangle$ with respect to the monomial ordering $<$ if any
composition $(f,g)_w$ of polynomials in $S$ is trivial modulo
$(S,w)$.

For a set $S\subseteq k\langle X\rangle$, the ideal of $k \langle
X\rangle$ generated by $S$ is denoted by $Id (S)$. If $S$ is a
Gr\"{o}bner-Shirshov basis in $k\langle X\rangle$ then it is also
called a Gr\"{o}bner-Shirshov basis for $Id(S)$ and for the algebra
$k\langle X\rangle/Id(S)=k\langle X|S\rangle$ generated by $X$ with
defining relations $S$.

The following lemma was first proved by Shirshov \cite{Sh} for free
Lie algebras  (with deg-lex ordering) (see also Bokut \cite{B72}).
Bokut \cite{B76} specialized the approach of Shirshov to associative
algebras  (see also Bergman \cite{Be}). For the case of commutative
polynomials, this lemma is known as the Buchberger's Theorem
\cite{Bu70}.

\begin{lemma}\label{CDL-alg}
 {\bf (Composition-Diamond Lemma)} \ Let $S$
be a subset of a free algebra $k\langle X\rangle$ over a field $k$,
and $<$ a monomial ordering on $X^*$. Then the following statements
are equivalent.
\begin{enumerate}
\item[ (i)] $S $ is a Gr\"{o}bner-Shirshov basis for $Id(S)$ with
respect to $<$.
\item[ (ii)] $f\in Id (S)\Rightarrow \bar{f}=a\bar{s}b$
for some $s\in S$ and $a,b\in  X^*$.
\item[ (iii)] $Irr (S) = \{ u \in X^* |  u \neq a\bar{s}b ,s\in S,a ,b \in X^*\}$
is a $k$-basis of the algebra $k\langle X | S \rangle.\ \square$
\end{enumerate}
\end{lemma}

If a subset $S$ of $k\langle X \rangle$ is not a
Gr\"{o}bner-Shirshov basis for $Id(S)$ then one can add to $S$ a
nontrivial composition $(f,g)_w $ of $f,g \in S$ and continue this
process repeatedly (actually using the transfinite induction) in
order to obtain a set $R$ of generators of $Id(S)$ such that any
composition of elements of $R$ is trivial modulo $R$ and the
corresponding ambiguity. Then $R$ is a Gr\"{o}bner-Shirshov basis
for $k\langle X|S\rangle$. This kind of process is called
\emph{Shirshov algorithm}.

Let $G=sgp\langle X|S\rangle$ be a semigroup presented by generators
$X$ and defining relations $S=\{u_i=v_i|i\in I\}$ for some index set
$I$. We will identify a semigroup relation $u=v\ ( u,v\in X^*)$ with
the algebra relation $u-v=0$ and with the polynomial (binomial)
$u-v\in k\langle X \rangle$. Then the semigroup algebra $kG$ has the
presentation $kG=k\langle X |S\rangle$. Because any composition of
binomials is again a binomial, from the Shirshov algorithm, it
follows that there exists a Gr\"{o}bner-Shirshov basis $R$ for $kG$
consisting of binomials such that $G=sgp\langle X |R\rangle$. Also,
$R$ does not depends on the field $k$. $R$ is also called a
Gr\"{o}bner-Shirshov basis for the semigroup $G$.

Actually, we do not need to use semigroup algebra $kG$ to define a
Gr\"{o}bner-Shirshov basis for a semigroup $G$. Let us reformulate
definition of a Gr\"{o}bner-Shirshov basis for a semigroup and
Composition-Diamond lemma for a free semigroup.

Let $<$  be  a monomial ordering on $X^*$, $S=\{(u,v)| v<u\}
\subseteq X^*\times X^*$  a set of semigroup relations, $u\equiv v \
mod(S)$ the congruence relation on $X^*$ generated by $S$,
$aub\mapsto avb$ and $avb\mapsto aub$ ( $(u,v)\in S,\ a,b\in X^*$)
the $S$-elementary transformations of $X$-words. Then $u\equiv v \ \
mod(S)$ if and only if there exists a sequence  $u=u_0\mapsto
u_1\mapsto \cdots \mapsto u_k=v$ of $S$-elementary transformations
of $u$ to $v$. Denote $S/_{\equiv \ mod(S)}=sgp\langle X|S\rangle$,
the semigroup generated by $X$ with defining relations $S$.

We will write $u\equiv v\ \  mod(S,w),\ w\in X^*$, if  $u=u_0\mapsto
u_1\mapsto \cdots \mapsto u_k=v$ for some $S$-elementary
transformations and $u_i<w\ (0\leq i\leq k)$.

Let $(u,v),\ (u',v')\in X^*\times X^*$ with $v<u,\ v'<u' $. We
define two kinds of compositions of $(u,v)\mbox{ and } (u',v')$:

$ (i)$ If  $w=ub=au'$ and $ |u|+|v|>|w|$ for some $a,b, w\in X^*$,
then $((u,v),\ (u',v'))_w=(vb,av')$ is called the \emph{intersection
composition} of $(u,v)\mbox{ and } (u',v')$ with respect to $w$.

$ (ii)$ If  $w=u=au'b$ for some $a,b,w\in X^*$, then $((u,v),\
(u',v'))_w=(v,av'b)$ is called the \emph{inclusion composition} of
$(u,v)\mbox{ and } (u',v')$ with respect to $w$.

A composition $(p,q)= ((u,v), (u', v'))_w$ is called \emph{trivial}
$mod(S,w)$, if $p\equiv q\ \  mod(S,w)$. It agrees with the above
definition of triviality of a composition for polynomials.

$S\subseteq X^*\times X^*$ is called a \emph{Gr\"{o}bner-Shirshov
basis} in $X^*$ relative to monomial ordering $<$ if any composition
of relations from $S$ is trivial
$mod(S,\textcolor[rgb]{0.0,0.00,0.00}{w})$. In this case, we also
call $S$ to be a Gr\"{o}bner-Shirshov basis for semigroup
$sgp\langle X|S\rangle$.

Then by Lemma \ref{CDL-alg}, we have
\begin{lemma}\label{CDL} {\bf (Composition-Diamond Lemma for semigroups)}
Let $X^*$ be a free monoid on set $X$ with a monomial
ordering $<$, $S=\{(u,v)| v<u\} \subseteq X^*\times X^*$. Then the
following conditions are equivalent.

($i$) \textcolor[rgb]{0.0,0.00,0.00}{S} is a Gr\"{o}bner-Shirshov
basis for the semigroup $sgp\langle X|S\rangle$ relative to $<$.

($ii$) If $p,q\in X^*,\ q<p \mbox{ and } p\equiv q \ mod (S)$, then
$p=aub$, where $(u,v)\in S$ and $v<u$ for some word $v$.

($iii$) $Irr(S) = \{p\in X^*|p\neq aub,(u,v)\in S \mbox{ for some
words}\  v, a,b\in X^*\}$ is a set of normal forms for $sgp\langle
X|S\rangle$.
\end{lemma}

\section{Ordered canonical idempotents}
We recall that an {\it inverse semigroup} is a semigroup in which
every element $a$ has a uniquely determined $a^{-1}$ such that
$aa^{-1}a=a\mbox{ and } a^{-1}aa^{-1}=a^{-1}$. The category of all
inverse semigroups (actually, it is a variety relative to two
operations $a\cdot b$ and $a^{-1}$) possesses free objects, free
inverse semigroups. Let $\mathcal {FI}(X)$ be a free inverse
semigroup generated by a set $X$, $X^{-1}=\{x^{-1}|x\in X\}$ with
$X\cap X^{-1}=\varnothing$. Denote $X\cup X^{-1}$ by $Y$. Then
$\mathcal {FI}(X)$ has the following semigroup presentation
$$
\mathcal {FI}(X)=sgp\langle Y|\ aa^{-1}a=a,\
aa^{-1}bb^{-1}=bb^{-1}aa^{-1}, \ a,b\in Y^* \rangle
$$
where $1^{-1}=1,\ (x^{-1})^{-1}=x\ (x\in X)$ and $(y_1y_2\cdots
y_n)^{-1}=y_n^{-1}\cdots y_2^{-1}y_1^{-1}\ (y_1,\ y_2,\cdots,\
y_n\in Y)$ (see, for example, \cite{how}, \cite{Pe}).

Let us assume that the set $Y$ is well ordered by an ordering $<$.
Let $<$ be also the corresponding deg-lex ordering of $Y^*$. For any
$u=y_1y_2\cdots y_n \ (y_1,\ y_2,\cdots,\ y_n\in Y)$, let
$fir(u)=y_1$.

We will define the {\it formal idempotents} in $Y^*$ which are
indeed the idempotents in the free inverse semigroup $\mathcal
{FI}(X)$. For the sake of convenience, we simply call the ``formal
idempotents" to be idempotents.

We give inductively definitions in $Y^*$ of an {\it idempotent,
canonical idempotent, prime canonical idempotent, ordered  (prime)
canonical idempotent} and {\it factors} of a canonical idempotent,
all of which but (prime) idempotent and ordered (prime) canonical
idempotent are defined in \cite{Schein}.

 (i) The empty word 1 is an idempotent, a canonical idempotent,
and an ordered canonical idempotent. This canonical idempotent has
no factors.

 (ii) If $h$ is an idempotent and $x\in Y$, then $x^{-1}hx$ is both an
idempotent and a prime idempotent. If $h$ is a canonical idempotent,
$x\in Y$ and the first letters of factors of $h$ are different from
$x$, then $x^{-1}hx$ is both a canonical idempotent and a prime
canonical idempotent. This canonical idempotent is its own factor.
Moreover, if the subword $h$ in this canonical idempotent is an
ordered canonical idempotent, then $x^{-1}hx$ is both an ordered
canonical idempotent and an ordered prime canonical idempotent.

 (iii) If $e_1, e_2,\cdots, e_m \ (m> 1)$ are prime idempotents,
then $e=e_1e_2\cdots e_m$ is an idempotent. Moreover, if $e_1,
e_2,\cdots, e_m$ are prime canonical idempotents and their first
letters are pairwise distinct, then $e=e_1e_2\cdots e_m$ is a
canonical idempotent and $e_1, e_2,\cdots, e_m$ are factors of $e$.
For this canonical idempotent, if $e_1, e_2,\cdots, e_m$ are ordered
canonical idempotents and $e\leq e_{i_1}e_{i_2}\cdots e_{i_m}$ for
any permutation $ (i_1, i_2,\cdots, i_m)$  of $ (1, 2,\cdots,m)$,
then $e$ is an ordered canonical idempotent.

\begin{remark} \label{permutation} By definition, it is easy to see that every
idempotent has even length. If $e=e_1e_2\cdots e_m$ is a canonical
idempotent, then $e$ is ordered if and only if
$fir(e_1)<fir(e_2)<\cdots<fir(e_m)$.
\end{remark}

\begin{lemma} (\cite{Schein})\label{Schein's}
Let $e=e_1e_2\cdots e_n\  (n\geq 1)$ be a canonical idempotent with
factors $e_1, e_2,\cdots, e_n$ and let $e_t=uv$ for some t  ($1\leq
t\leq n$) with $u, v \neq 1$. Then neither $e_1\cdots e_{t-1}u$ nor
$ve_{t+1}\cdots e_n$ is a canonical idempotent.
\end{lemma}

The following lemma is a generalization of Lemma \ref{Schein's}.

\begin{lemma}\label{unique}
Let $e=e_1e_2\cdots e_i\cdots e_j\cdots e_n \ (1\leq i<j\leq n)$ be
a canonical idempotent with factors $e_1, e_2,\cdots, e_n$ and let
$e_i=u_iv_i, \ e_j=u_jv_j$ with either $u_i,v_i\neq1$ or
$v_j,u_j\neq1$. Then
 $w=v_ie_{i+1}\cdots e_{j-1}u_j$ is not a canonical
idempotent.
\end{lemma}
\begin{proof}
We prove the lemma by induction on $k=|w|$, the length of $w$.

If $k=1$, then $w=x\in Y$ and the result holds. Suppose that the
result holds for all $w$ with $|w|\leq l$. Consider $w$ with
$|w|=l+1$.

If one of $u_i, v_i, v_j, u_j$ is empty, then our statement holds by
Lemma \ref{Schein's}. Now we suppose that $u_i, v_i, v_j, u_j\neq
1$. By way of contradiction, assume that $w=w_1w_2\cdots w_s$ is a
canonical idempotent with factors $w_1, w_2,\cdots, w_s$. By Lemma
\ref{Schein's}, $v_i$ is not a canonical idempotent, and hence
$v_i=w_1\cdots w_{k-1}a_k \ (1\leq k\leq s)$, where $w_k=a_kc_k$ and
$a_k, c_k\neq 1$. Similarly, $u_j=b_lw_{l+1}\cdots w_s \ (k\leq
l\leq s)$ for $w_l=d_lb_l$ and $d_l, b_l\neq 1$.

{\it Case 1}. $k=l$. Then $w_k=w_l=x^{-1}hx=a_ke_{i+1}\cdots
e_{j-1}b_l=x^{-1}a_k' e_{i+1}\cdots e_{j-1}b_l'x$, where $x\in Y,
a_k=x^{-1}a_k', b_l=b_l'x $ and $h$ is a canonical idempotent. If
$a_k'b_l'=1$, then $fir(e_i)=fir(e_j)=x$, which is impossible since
$e$ is a canonical idempotent. Thus $a_k'b_l'\neq1$ and by induction
hypothesis, $h=a_k' e_{i+1}\cdots e_{j-1}b_l'$ is not a canonical
idempotent since $|h|<|w|$, which is a contradiction.

{\it Case 2}. $k<l$. Then, by induction hypothesis, $e_{i+1}\cdots
e_{j-1}=c_kw_{k+1}\cdots w_{l-1}d_l$ is not a canonical idempotent
since $a_k, c_k\neq1$, which is also a contradiction.
\end{proof}

\begin{lemma}\label{subword}
Suppose that $w\in Y^*$ is an idempotent. Then $w$ is a canonical
idempotent if and only if $w$ has no subword of the form
$x^{-1}exfx^{-1}$, where $x\in Y$,  $x^{-1}ex$ and $xfx^{-1}$ are
both prime canonical idempotents.
\end{lemma}
\begin{proof}
We use induction on $k$, where $2k=|w|$.

We first prove the ``if" part. If $k=0$, then the ``if" part clearly
holds. Suppose that the ``if" part holds for all $w$ with $|w|\leq
2l$. Consider $w$ with $|w|=2l+2$. Suppose that $w=w_1w_2\cdots
w_s$, where $w_1,w_2,\cdots, w_s$ are prime idempotents. If $s=1$,
then $w=y^{-1}hy$, where $y\in Y$ and $h$ is an idempotent. By
induction hypothesis, $h$ is a canonical idempotent. Since $w$ has
no subwords of the form $x^{-1}exfx^{-1}$, the first letters of the
factors of $h$ are not $y$. Thus $w$ is a canonical idempotent. If
$s\geq 2$, then by induction hypothesis, $w_1,w_2,\cdots ,w_s$ are
all canonical idempotents, and the first letters of $w_1,w_2,\cdots
,w_s$ are pairwise distinct. Hence $w$ is a canonical idempotent.

Now we prove the ``only if" part. If $k=0$, then the ``only if" part
holds. Suppose that the ``only if" part holds for all $w$ with
$|w|\leq 2l$. Consider $w$ with $|w|=2l+2$. By way of contradiction,
we assume that $w=w_1w_2\cdots w_s$ with factors $w_1,w_2,\cdots,
w_s$ and subword $x^{-1}exfx^{-1}$, where $x^{-1}ex$ and $xfx^{-1}$
are both prime canonical idempotents.

If $s=1$, then $w=y^{-1}hy=y^{-1}h_1h_2\cdots h_ky$, where $y\in Y$
and $h$ is a canonical idempotent with factors $h_1,h_2,\cdots,h_k \
(k\geq1)$. By induction hypothesis, $x^{-1}exfx^{-1}$ is not subword
of $h$, and then $x^{-1}exfx^{-1}$ is a beginning or end part of
$w$. For the former case,  by Lemma \ref{unique}, we have that
$x=y$, $e=h_1\cdots h_i$ and $xfx^{-1}=h_{i+1}$ for some $i$, and
hence $fir(h_{i+1})=x=y$, which is a contradiction since $w$ is a
canonical idempotent. Similarly, we can get a contradiction for the
latter case.

If $s>1$, then, by induction hypothesis, $x^{-1}exfx^{-1}$ is not a
subword of $w_2\cdots w_s$ or $w_1\cdots w_{s-1}$ and so
$x^{-1}exfx^{-1}=x^{-1}v_1w_2\cdots w_{s-1}u_sx^{-1}$,
$w_1=u_1x^{-1}v_1$ and $w_s=u_sx^{-1}v_s$ for some
$u_1,v_1,u_s,v_s\in Y^*$. Since $w$ is canonical, $v_1v_s\neq1$.
Then, by Lemma \ref{unique}, $exfx^{-1}=v_1w_2\cdots
w_{s-1}u_sx^{-1}$ is not a canonical idempotent, a contradiction.
\end{proof}

\begin{lemma}\label{insert}
 Let $e'$ be an idempotent, $a,b\in Y^*$. Then $e=ab$ is an
idempotent if and only if $ae'b$ is an idempotent.
\end{lemma}
\begin{proof} We may assume that $e'$ is a nonempty idempotent.

We first prove the ``if'' part. Ordering the set $\{(a,b)|a,b\in
\mathbb{Z}_+\}$ lexicographically, we prove the ``if'' part by
induction on $(|ae'b|,|e'|)$. If $(|ae'b|,|e'|)=(2,2)$, then $ab=1$
is an idempotent. Suppose that the ``if'' part holds for all
$a,b,e'$ with $(|ae'b|,|e'|)<(2l,2k),l,k\geq1$. Consider $a,b,e'$
with $(|ae'b|,|e'|)=(2l,2k)$ and $ab\neq1$. Suppose that
$ae'b=e_1e_2\cdots e_m (m\geq 1)$, where $e_1, e_2,\cdots, e_m$ are
prime idempotents.

{\it Case 1.} $|e'|>2$, i.e., $e'=ce''d$ with some nonempty
idempotent $e''$ as a proper subword. Then, by induction hypothesis,
$acdb$ and $cd$ are idempotents and so is $ab$.

{\it Case 2.} $|e'|=2$, i.e., $e'=xx^{-1},\  x \in Y$.

{\it Subcase 1.} $e_i=ce'd$  for some $c,d\in Y^*$ and $1\leq i\leq
m$. If $m=1$, then we may suppose $ae'b=yf_1\cdots f_py^{-1}\
(p\geq1)$, where $f_1,\cdots, f_p$ are prime idempotents. Moreover,
if $a=1$ ($b=1$ is similar), i.e., $x=y,\ ae'b=xx^{-1}f_1'xf_2\cdots
f_px^{-1}$, where $f_1'$ is an idempotent, then $ab=f_1'xf_2\cdots
f_px^{-1}$ is an idempotent. If $a\neq1$ and $b\neq1$, then
$f_1\cdots f_p=ce'd$ for some $c,d\in Y^*$. Hence, $cd$ and
$ab=ycdy^{-1}$ are both idempotents by induction hypothesis and by
definition respectively.

If $m>1$, then $cd$ and $ab=e_1\cdots e_{i-1}cde_{i+1}\cdots e_m$
are both idempotents by induction hypothesis and by definition
respectively.

{\it Subcase 2.} $e_i=x^{-1}e_i'x$ and $e_{i+1}=x^{-1}e_{i+1}'x$
($1\leq i\leq m-1$), where $e_i'$ and $e_{i+1}'$ are idempotents,
i.e., $e_ie_{i+1}=x^{-1}e_i'e'e_{i+1}'x$. Then, $ab=e_1\cdots
e_{i-1}x^{-1}e_i'e_{i+1}'xe_{i+2}\cdots e_m$ is an idempotent by
definition.

Now we prove the ``only'' part. We also prove it by induction on
$k$, where $2k=|e|$. If $k=0$, then $ae'b=e'$ is an idempotent.
Suppose that the ``only'' part holds for all $e$ with $|e|\leq 2l$.
Consider $e$ with $|e|=2l+2$. Suppose that $e=e_1e_2\cdots e_m\
(m\geq 1)$, where $e_1, e_2,\cdots, e_m$ are prime idempotents.

{\it Case 1.}  $m=1$. Then $e=x^{-1}hx$, where $x\in Y$ and $h$ is
an idempotent. If $a=1$ or $b=1$, then our statement holds by
definition. If $a,b\neq1$, then we suppose that $a=x^{-1}a'$ and
$b=b'x$. Now, by induction hypothesis, $a'e'b'$ is clearly an
idempotent, and  $ae'b=x^{-1}a'e'b'x$ is also an idempotent.

{\it Case 2.} $m>1$. Then $ae'b=e_1\cdots e_{i-1}ce'de_{i+1}\cdots
e_m\ (1\leq i\leq m)$, where $cd=e_i$. By induction hypothesis,
$ce'd$ is an idempotent and so is $ae'b$.
\end{proof}

\begin{lemma}\label{ordered}
If $w,e$ and $f$ are nonempty ordered canonical idempotents and
$w=aefb$ for some $a,b\in Y^*$, then $ef<fe$.
\end{lemma}
\begin{proof}
We prove the lemma by induction on $k$, where $2k=|w|\geq4$. If
$k=2$, then $w=ef$ and our statement holds.
 Suppose that this lemma holds for
all $w$ with $|w|\leq 2l$. Consider $w$ with $|w|=2l+2$. Suppose
that $w=w_1w_2\cdots w_m$ with factors $w_1, w_2, \cdots, w_m$. If
$m=1$, then $w=x^{-1}hx=x^{-1}h_1\cdots h_nx$, where $x\in Y$ and
$h$ is a canonical idempotent with factors $h_1, \cdots, h_n$. By
Lemma \ref{unique}, $ef$ is a subword of $h$, and by induction
hypothesis $ef<fe$. If $m>1$, then by Lemma \ref{unique}, $ef$ is
either a product of factors of $w$ or a subword of some $w_i$ for
$1\leq i\leq m$. Hence, we have $ef<fe$ by definition in the former
case or by induction hypothesis in the latter case.
\end{proof}

\ \

In Lemmas \ref{algo} and \ref{equivalent} as well as in Proposition
\ref{pro},  $S\subset Y^*\times Y^*$ denotes the set of the
following defining relations (a) and (b):

 (a) $(ef,fe)$, where both $e$ and $f$ are ordered prime canonical idempotents,
$ef$ is a canonical idempotent and $fe<ef$;

(b) $(x^{-1}e'xf'x^{-1}, f'x^{-1}e')$, where $x\in Y$, both
$x^{-1}e'x$ and $xf'x^{-1}$ are ordered prime canonical idempotents.

\begin{lemma}\label{algo}
(1) Suppose that $e$ is a prime canonical idempotent, $w\in Y^*$ and
$e<w$. Then, there exists a prime ordered canonical idempotent $e'$
such that $fir(e)=fir(e')$ and $e\equiv e'\ mod(S,w)$.

(2) Suppose that $e=e_1e_2\cdots e_m\ (m\geq1)$ is a canonical
idempotent with factors $e_1,e_2,\cdots, e_m$, $w\in Y^*$ and $e<w$.
Then, there exists an ordered canonical idempotent
$e'=e_{i_1}e_{i_2}\cdots e_{i_m}$ such that $e\equiv e'\ mod(S,w)$,
where $(i_1, i_2,\cdots, i_m)$ is a permutation of $(1,
2,\cdots,m)$.

(3) Suppose that $e$ is a nonempty idempotent, $w\in Y^*$ and $e<w$.
Then, there exists a canonical idempotent $e'$ such that $e\equiv
e'\ mod(S,w)$.
\end{lemma}
\begin{proof}
(1). We use induction on $k=|e|$. If $k=2$, then, by taking $e'=e$,
(1) holds. Suppose that (1) holds for all prime canonical idempotent
$e$ with $|e|\leq 2l$. We consider $e$ with length $2l+2$. Now, by
induction hypothesis, we may suppose $e=x^{-1}e_1e_2\cdots e_mx$,
where $x\in Y$ and $e_1,e_2,\cdots ,e_m$ are ordered prime canonical
idempotents. If $e$ is ordered, then (1) holds. Assume $e$ is not
ordered, i.e., $e_1e_2\cdots e_m$ is not ordered (so $m>1$). By
Remark \ref{permutation}, there exists a permutation $(i_1,
i_2,\cdots, i_m)$ of $(1, 2,\cdots,m)$ such that
$e'=x^{-1}e_{i_1}e_{i_2}\cdots e_{i_m}x$ is an ordered canonical
idempotent.

It suffices to prove that $e_1e_2\cdots e_m\equiv
e_{i_1}e_{i_2}\cdots e_{i_m}\ mod(S,w')$ for any word $w'$ such that
$e_1e_2\cdots e_m<w'$. We prove it by induction on $m$. If $m=2$,
then our statement holds clearly. Supposing our statement holds for
$m\leq n$, we consider $m=n+1$. If $e_1\neq e_{i_1}$, i.e., $1<i_1$,
then $fir(e _{i_1})<fir(e _t)$ for $1\leq t<i_1$. Hence,
$mod(S,w')$, the following $\equiv$'s hold,
\begin{eqnarray*}
e_1 e_2 \cdots e_m &\equiv& e_1e_2\cdots
e_{i_1-3} e_{i_1-2} e_{i_1} e_{i_1-1} e_{i_1+1} \cdots e_m \\
&\equiv&e_1e_2\cdots
e_{i_1-3} e_{i_1} e_{i_1-2} e_{i_1-1} e_{i_1+1} \cdots e_m \\
&\vdots&\\
&\equiv&e_{i_1}e_1 e_2\cdots e_{i_1-3} e_{i_1-2} e_{i_1-1} e_{i_1+1}
\cdots e_m.
\end{eqnarray*}
Thus, we may suppose $e_1=e_{i_1}$. Then, by induction hypothesis,
$e_1e_2\cdots e_m\equiv e_{i_1}e_{i_2}\cdots e_{i_m}$.

This ends our proof of (1).

(2) follows from the proof of (1).

(3). We use induction on $|e|$. If $|e|=2$, then, by taking
$e'=e=x^{-1}x$, (1) holds. Suppose that (1) holds for all idempotent
$e$ with $|e|\leq 2l$. We consider $e$ with length $2l+2$. If $e$ is
not canonical, then, by Lemma \ref{subword}, $e=ay^{-1}fygy^{-1}b$,
where  $y\in Y,\ a,b\in Y^*$, $y^{-1}fy$ and $ygy^{-1}$ are both
canonical idempotents. By (1), we may suppose $y^{-1}fy$ and
$ygy^{-1}$ are both ordered canonical idempotents. Then, $e\equiv
agy^{-1}fb\equiv e'\ mod(S,w)$, where $e'$ is a canonical
idempotent, and the second $\equiv$ holds by induction hypothesis
since $agy^{-1}fb$ is an idempotent by Lemma \ref{insert}.
\end{proof}

\begin{lemma}\label{equivalent}

 (1) Suppose that $e$ and $f$ are both idempotents and
$x^{-1}exfx^{-1}<w$ for some $x \in Y,\ w\in Y^*$. Then
$x^{-1}exfx^{-1}\equiv fx^{-1}e\ \ mod (S,w)$.

 (2) Suppose that $e$ and
$f$ are both nonempty idempotents and $ef,fe<w$
 for some $w\in Y^*$. Then $ef\equiv fe\ \
mod (S,w)$.
\end{lemma}
\begin{proof}
(1). We use induction on $k=|x^{-1}exfx^{-1}|\geq3$. If $k=3$, then
$e=f=1$ and (1) holds. Supposing (1) holds for all $x^{-1}exfx^{-1}$
with $k\leq2l-1$, we consider $x^{-1}exfx^{-1}$ with $k=2l+1$. By
Lemma \ref{algo}, we may suppose $e$ and $f$ are both ordered
canonical idempotents. If $x^{-1}ex$ and $xfx^{-1}$ are both
canonical, then (1) holds. If $x^{-1}ex$ or $xfx^{-1}$ is not
canonical, say, $x^{-1}ex$ is not canonical, then $e=e_1\cdots
e_{i-1}xgx^{-1}e_{i+1}\cdots e_n\ (n\geq1)$ for some integer $i$,
where $e_1,\cdots, e_{i-1},\ xgx^{-1}=e_i,\ e_{i+1},\cdots ,e_n$ are
factors of $e$. Hence, $mod(S,w)$, the following $\equiv$'s hold by
induction hypothesis,
\begin{eqnarray*}
x^{-1}exfx^{-1}&=&x^{-1}e_1\cdots e_{i-1}xgx^{-1}e_{i+1}\cdots
e_nxfx^{-1}\\
&\equiv& gx^{-1}e_1\cdots e_{i-1}e_{i+1}\cdots
e_nxfx^{-1}\\
&\equiv& gfx^{-1}e_1\cdots e_{i-1}e_{i+1}\cdots e_n.
\end{eqnarray*}
On the other hand, we have
\begin{eqnarray*}
fx^{-1}e&=&fx^{-1}e_1\cdots e_{i-1}xgx^{-1}e_{i+1}\cdots e_n\\
&\equiv&fgx^{-1}e_1\cdots e_{i-1}e_{i+1}\cdots e_n.
\end{eqnarray*}
Now, it suffices to prove that $gf\equiv fg\ \ mod(S,w')$ for any
word $w'$ such that $max\{gf,fg\}<w$. We prove it by induction on
$t=|gf|$. Suppose $g=g_1g_2\cdots g_m$ and $f=g_{m+1}g_{m+2}\cdots
g_{m+n}$, where $m,n\geq1$ and each $g_j\ (1\leq j\leq m+n)$ is
prime ordered canonical idempotent. If $t=0$, then $gf\equiv fg\
mod(S,w')$. Supposing $gf\equiv fg\ mod(S,w')$ for any $g$ and $f$
with $t\leq2l$, we consider $g$ and $f$ for $t=2l+2$. If $gf$ is
canonical, then by Lemma \ref{algo}, there exists a permutation
$(i_1,i_2,\cdots, i_{m+n})$ of $(1,2,\cdots, m+n)$ such that
$gf\equiv g_{i_1}\cdots g_{i_{m+n}}\equiv fg\ mod(S,w')$. If $gf$ is
not canonical, i.e., $g_s=x^{-1}g_s'x$ and $g_{m+j}=x^{-1}g_{m+j}'x$
for some $x\in Y$, integers $s,j \ (1\leq s\leq m,\ 1\leq j\leq n)$,
and ordered canonical idempotents $g_s'$ and $g_{m+j}'$, then,
$mod(S,w')$, the following $\equiv$'s hold by induction hypothesis
on $|gf|$ or on $|x^{-1}exfx^{-1}|$,
\begin{eqnarray*}
 gf&=&g_1\cdots g_{s-1}x^{-1}g_s'xg_{s+1}\cdots g_mg_{m+1}\cdots
 g_{m+j-1}x^{-1}g_{m+j}'xg_{m+j+1}\cdots g_{m+n}\\
 &\equiv&g_1\cdots g_{s-1}g_{s+1}\cdots g_mg_{m+1}\cdots
 g_{m+j-1}x^{-1}g_s'g_{m+j}'xg_{m+j+1}\cdots g_{m+n}\\
 &\equiv&g_{m+1}\cdots
 g_{m+j-1}g_{m+j+1}\cdots g_{m+n}g_1\cdots g_{s-1}x^{-1}g_{m+j}'g_s'xg_{s+1}\cdots g_m\\
 &\equiv&g_{m+1}\cdots
 g_{m+j-1}x^{-1}g_{m+j}'xg_{m+j+1}\cdots g_{m+n}g_1\cdots g_{s-1}x^{-1}g_s'xg_{s+1}\cdots g_m\\
 &=&fg.
\end{eqnarray*}

This shows (1).

 (2). By Lemma \ref{algo}, we may
assume that $e$ and $f$ are both ordered canonical idempotents.
Then, (2) follows from the proof of (1).
\end{proof}

\ \

By Lemma \ref{equivalent}, for any $a,b\in Y^*$, $aa^{-1}a=a$ and
$aa^{-1}bb^{-1}=bb^{-1}aa^{-1}$ in $sgp\langle Y|S\rangle$. Thus, we
have the following proposition.

\begin{proposition}\label{pro}
$sgp\langle Y|S\rangle$ is a free inverse semigroup with identity,
i.e., $sgp\langle Y|S\rangle=\mathcal {FI}(X)$.
\end{proposition}

\section{Main theorem}

The following theorem is the main result of this paper.

\begin{theorem}\label{gsb-sgp}
Let $X$ be a set, $X^{-1}=\{x^{-1}|x\in X\}$ with $X\cap
X^{-1}=\varnothing$, $<$ a well ordering on $Y=X\cup X^{-1}$ and
also the deg-lex ordering of $Y^*$, $\mathcal {FI}(X)$ the free
inverse semigroup generated by $X$. Let $S\subset Y^*\times Y^*$ be
the set of the following defining relations (a) and (b):

 (a) $(ef,fe)$, where both $e$ and $f$ are ordered prime canonical idempotents,
$ef$ is a canonical idempotent and $fe<ef$;

(b) $(x^{-1}e'xf'x^{-1}, f'x^{-1}e')$, where $x\in Y$, both
$x^{-1}e'x$ and $xf'x^{-1}$ are ordered prime canonical idempotents.

Then $S$ is a Gr\"{o}bner-Shirshov basis for the free inverse
semigroup $\mathcal {FI}(X)$ with respect to the deg-lex ordering of
$Y^*$.
\end{theorem}
\begin{proof}
In the following, all $\equiv$'s hold $mod(S,w)$ by Lemma
\ref{insert} or/and Lemma \ref{equivalent}. By notation (a$\wedge$b)
we denote all the possible compositions between the relations of
type (a) and of type (b) in $S$, and similarly we use notations
(a$\wedge$a), (b$\wedge$a) and (b$\wedge$b).

We check all the possible compositions step by step. For
convenience, we use the algebra language of Lemma \ref{CDL-alg}
rather than the semigroup language of Lemma \ref{CDL}.

 \ \

 (a$\wedge$a)\ \ \ \  $ef-fe\wedge e'f'-f'e'$.

 (1) Inclusion compositions.

 By Lemma \ref{ordered}, $e'f'$ can not be a subword of $e$ or
 $f$. Thus, by Lemma \ref{unique}, there are
 no inclusion compositions.

 (2) Intersection compositions.

 There are five cases to
consider.

{\it Case 1.} $e=ae'b$ for some $a,b\in Y^*$. Then $w=ae'bfc$ and
\begin{eqnarray*}
 (ef,e'f')_w
&=&-fae'bc+abfce'\\
&\equiv&-fabce'+fabce'\\
&\equiv&0
\end{eqnarray*}

{\it Case 2.} $e=ab,\ f=cd,\ e'=bc$ for some $a,b,c,d\in Y^*$ and
$b\neq1$. By Lemma \ref{unique}, this case is impossible.

In the following cases, similar to the Case 1, $(ef,e'f')_w\equiv0\
\ mod (S,w)$. We list only the ambiguity $w$ for each case.

{\it Case 3.} $e=ab,\ e'=bfc$ for some $a,b,c\in Y^*$. Then
$w=abfcf'$.

{\it Case 4.} $f=ae'b,\ f'=bc$ for some $a,b,c\in Y^*$. Then
$w=eae'bc$.

{\it Case 5.} $f=ab,\ e'=bc$ for some $a,b,c\in Y^*$. Then
$w=eabcf'$.

\ \

 (a$\wedge$b)\ \ \ \  $ef-fe\wedge x^{-1}e'xf'x^{-1}-f'x^{-1}e'$.

By Lemma \ref{subword}, there are no
 inclusion compositions of type a$\wedge$b. To consider the
intersection compositions, there are five cases to consider.

{\it Case 1}. $e=ax^{-1}e'xb,\ f'=bfc$ for some $a,b,c\in Y^*$. Then
$w=ax^{-1}e'xbfcx^{-1}$ and
\begin{eqnarray*}
(ef,x^{-1}e'xf')_w
&=&-fax^{-1}e'xbcx^{-1}+abfcx^{-1}e'\\
&\equiv&-fabcx^{-1}e'+fabcx^{-1}e'\\
&\equiv&0
\end{eqnarray*}

{\it Case 2}. $e=ax^{-1}b,\ f=cxd,\ e'=bc,\ f'=dg$ for some
$a,b,c,d,g\in Y^*$ and $b\neq1$. By Lemma \ref{unique}, this case is
impossible.

In the following cases, similar to the Case 1,
$(ef,x^{-1}e'xf'x^{-1})_w\equiv0\ \ mod (S,w)$. We list only the
ambiguity $w$ for each case.

{\it Case 3}. $e=ax^{-1}b,\ e'=bfc$ for some $a,b,c\in Y^*$. Then
$w=ax^{-1}bfcxf'x^{-1}$.

{\it Case 4}. $f=ax^{-1}e'xb,\ f'=bc$ for some $a,b,c\in Y^*$. Then
$w=eax^{-1}e'xbcx^{-1}$.

{\it Case 5}. $f=ax^{-1}b,\ e'=bc$ for some $a,b\in Y^*$. Then
$w=eax^{-1}bcxf'x^{-1}$.

\ \

 (b$\wedge$a)\ \ \ \  $x^{-1}exfx^{-1}-fx^{-1}e\wedge e'f'-f'e'$.

 (1) Inclusion compositions.

By Lemma \ref{ordered}, $e'f'$ can not be a subword of $x^{-1}ex$ or
$xfx^{-1}$. Hence, $e'f'$ is a subword of $x^{-1}exf$ or
$exfx^{-1}$, and by Lemma \ref{unique}, $e'f'=x^{-1}exf_1$ or
$e'f'=e_mxfx^{-1}$, where $e=e_1e_2\cdots e_m$ and $f=f_1f_2\cdots
e_n$ with factors $e_1, e_2, \cdots,  e_m$ and $f_1, f_2, \cdots,
f_n$ respectively. Now, it is easy to check that all the inclusion
compositions are trivial.

(2) Intersection compositions.

This case is symmetrical to the case of intersection compositions of
type a$\wedge$b.

\ \

 (b$\wedge$b)\ \ \ \  $x^{-1}exfx^{-1}-fx^{-1}e\wedge y^{-1}e'yf'y^{-1}-f'y^{-1}e'$.

 (1) Inclusion compositions.

 By Lemma \ref{subword},
$y^{-1}e'yf'y^{-1}$ can not be a subword of $x^{-1}ex$ or
$xfx^{-1}$. Then, there are on inclusion compositions.

 (2) Intersection compositions.

 There are six cases to consider.

{\it Case 1.} $e=ay^{-1}e'yb, f'=bxfx^{-1}c$ for some $a,b,c\in
Y^*$. Then, $ w=x^{-1}ay^{-1}e'ybxfx^{-1}cy^{-1} $ and
\begin{eqnarray*}
 (x^{-1}exfx^{-1},y^{-1}e'yf'y^{-1})_w&=&-fx^{-1}ay^{-1}e'ybcy^{-1}+x^{-1}abxfx^{-1}cy^{-1}e'\\
&\equiv&-fx^{-1}abcy^{-1}e'+fx^{-1}abcy^{-1}e'\\
&\equiv&0
\end{eqnarray*}

In the following cases, we also have
$(x^{-1}exfx^{-1},y^{-1}e'yf'y^{-1})_w\equiv 0$ in a similar way.

{\it Case 2.} $e=ab, y^{-1}e'y=bx^{-1}c,f'=dx^{-1}g$ for some
$a,b,c,d\in Y^*$. Since $y^{-1}e'y$ is a prime canonical idempotent,
by Lemma \ref{unique}, $b=d=1$, i.e., $x=y^{-1}$ and $f'=e'$. Thus,
$w=x^{-1}exfx^{-1}f'x$.

{\it Case 3.} $e=ab, y^{-1}e'=bx^{-1}fc$ for some $a,b,c\in Y^*$. By
Lemma \ref{unique}, $b\neq1$. Suppose $b=y^{-1}b'$ for some $b'\in
Y^*$. Then we have $w=x^{-1}ay^{-1}b'xfx^{-1}cyf'y^{-1}$.

{\it Case 4.} $f=ay^{-1}e'yb, f'=bx^{-1}c$ for some $a,b,c\in Y^*$.
Then, $w=x^{-1}exay^{-1}e'ybx^{-1}cy^{-1}$.

{\it Case 5.} $f=ay^{-1}b, e'=bx^{-1}c$ for some $a,b,c\in Y^*$.
Then, $w=x^{-1}exay^{-1}bx^{-1}cyf'y^{-1}$.

{\it Case 6.} The intersection of $x^{-1}exfx^{-1}$ and
$y^{-1}e'yf'y^{-1}$ is $x^{-1}=y^{-1}$. Then,
$w=x^{-1}exfx^{-1}e'xf'x^{-1}$.

 Therefore, all possible
compositions in $S$ are trivial.
\end{proof}

 \ \

By Theorem \ref{gsb-sgp} and Composition-Diamond Lemma, $Irr (S)$ is
normal forms of the free inverse semigroup $\mathcal {FI}(X)$. It is
easy to see that $Irr (S) = \{ u \in (X\cup X^{-1})^* | u \neq
a\bar{s}b ,s\in S,a ,b \in (X\cup X^{-1})^*\}$ consists of the word
$u_0e_1u_1\cdots e_mu_m\in (X\cup X^{-1})^* $, where $m\geq0, \
u_1,\cdots, u_{m-1}\neq1$, $u_0u_1\cdots u_m$ has no subword of form
$yy^{-1}$ for $y\in X\cup X^{-1}$, $e_1,\cdots, e_m$ are ordered
canonical idempotents, and the first (last, respectively) letters of
the factors of $e_i\ (1\leq i\leq m)$ are not equal to the first
(last, respectively) letter of $u_i\ (u_{i-1},\
\mbox{respectively})$. Thus $Irr (S)$ is a set of canonical words in
the sense of \cite{Schein}, and different words in $Irr (S)$
represent different elements in $\mathcal {FI}(X)$.

\ \

\noindent{\bf Acknowledgement}: The authors would like to express
their deepest gratitude to L. N. Shevrin  for giving us some
valuable comments and useful remarks.

\ \

\end{document}